\RequirePackage{fix-cm}
\documentclass[twocolumn]{svjour3}          
\smartqed  
\usepackage{graphicx}
\usepackage{subcaption}
\usepackage{amsmath}
\usepackage{verbatim}
\usepackage{hyperref}
\usepackage{amsrefs}
\usepackage{float}
\usepackage{amssymb}
\usepackage{bbm}
\usepackage{bbold}
\usepackage{comment}
\usepackage{xcolor}
\usepackage{multicol}
\usepackage{neuralnetwork}
\usepackage{cancel}
\usepackage{array}
\usepackage{mathptmx}      
\usepackage{latexsym}
\DeclareMathOperator*{\argmin}{arg\,min}
\newcommand{\bXi}{\mbox{\boldmath $\Xi$}}
\newcommand{\bxi}{\mbox{\boldmath $\xi$}}
\newcommand{\bTheta}{\mbox{\boldmath $\Theta$}}
\journalname{Nonlinear Dynamics}
\begin{document}

\title{Compressed Compressor
}


\author{Alyssa Novelia \textsuperscript{1}       \and
        Yusuf Aydogdu \textsuperscript{1} \and
        Thambirajah Ravichandran \textsuperscript{1} \and
        N. Sri Namachchivaya \textsuperscript{1}
}


\institute{\textsuperscript{1} Corresponding author: N. Sri Namachchivaya \at
              Department of Applied Mathematics, University of Waterloo\\
              Tel.: +1 (519) 888-4567 Ext. 38156 \\
              \email{nsnamachchivaya@uwaterloo.ca}           
           }

\date{Received: date / Accepted: date}

\maketitle

\begin{abstract}

In this paper, we present a dada-driven reduced order model of viscous Moore-Greitzer~(MG) partial differential equation (PDE) by threading together ideas from principal component analysis (PCA) and autoencoder neural networks to sparse regression and compressed sensing. 
Numerical simulation of the infinite dimensional viscous MG system is reduced into low dimensional data using PCA and autoencoder neural networks based reduced order modelling (ROM) approaches. 
Based on the observation that MG equations close to bifurcations have a sparse representation (normal form) with respect to high-dimensional polynomial spaces, we use the Sparse Identification of Dynamical Systems (SINDy) algorithm which uses a collection of all monomials as sampling matrix and the LASSO algorithm to recover a system of sparse two ordinary differential equations (ODEs) with cubic nonlinearities.
The discovered governing equations can be used to fully recover the original system dynamics up to 98.9\% accuracy.
When dimensional reduction is performed along the dataset's principal components, the resulting low dimensional differential equations will be consistent and have some resemblance to the normal form structure.
Additionally,  a new nonlinear behaviour is exhibited in viscous MG equations during rotating stall instability past the Hopf bifurcation point.


\end{abstract}

\section{Introduction}
\label{intro}
This paper develops data-driven theory and algorithms to detect and mitigate stall compressor instability. The motivation is to produce a high-fidelity simulation of a jet engine compressor called the digital twin, which has the ability to monitor and diagnose complex systems to improve performance efficiency and utilization. Jet engine compressor models typically integrate a hierarchy of multi-physics and multi-fidelity models which are continually updated with data streams from the sensors. The model used to describe airflow inside the jet engine compressor is the viscous MG equations \cites{GM86-1,GM86-2} which consist of a nonlinear partial differential equations(PDE) \eqref{laplace-disturbance} and two ODEs \eqref{phi-eqn-viscous-avg} and \eqref{pressure-ode}. The PDE describes the spatiotemporal behavior of disturbances in the inlet region of the compression system and the two ODEs describe the coupling of the disturbances with the mean flow and pressure. There are three types of Hopf bifurcations that can exist in the viscous MG equations corresponding to physical oscillations dominated by the ODE (surge), PDE (rotating stall), or a mixture of both. The objective is to use optimization and regression techniques from machine learning to arrive at a lower dimensional description of the PDE from datasets - hence the name "compressed compressor". The success of compressed compressor is rooted on accurate representations of the multi-physics and multi-fidelity models.

In Section \ref{sec:s2} we introduce the viscous MG equations, provide an explicit expression for the system's equilibrium, and show that the steady operating axial flow and pressure drifts from the aforementioned equilibrium during PDE bifurcation. In Section \ref{sec:s3}, we introduce reduced-order modeling (ROM) to significantly alleviate computational costs by projecting the high dimensional state variables onto a low-dimensional subspace. We perform ROM on simulated data from viscous MG equations to construct a set of "good" basis functions. Approximations of bases spanning this subspace are constructed using principal component analysis (PCA) \cites{Sirovich87,CZN07, HLBR12} and both linear and nonlinear autoencoder neural networks \cites{Plaut18, Kramer91, Scholz02}.

It is impossible to effectively ``learn" from high dimensional data unless there is some kind of implicit or explicit low dimensional structure - for which there are multiple mathematically precise definitions. Over the past 10 years, researchers have focused on sparsity as one type of criteria for low-dimensional structure. The inherent sparsity of natural signals is central to the mathematical framework of compressed sensing \cites{Donoho06, CT05, CRT06}. The main aim of compressed sensing is to construct a sparse vector from linear measurements of the vector such that the number of observed measurements $m$ is significantly smaller than the dimension $n$ of the original vector and satisfies the ``Restricted Isometry Property" (RIP). Intuitively, the existence of a RIP implies that the geometry of sparse vectors is preserved through the measurement matrix. These techniques rely heavily on the fact that many dynamical systems can be represented by governing equations that are sparse in the space of all possible functions. The assumption for the low dimensional structure for the MG equations originates from the center manifold theory  in dynamical systems \cites{GH83, XB00}, where a high dimensional system undergoing Hopf bifurcation can be fully described by projecting the equations onto the subspace of a 2-dimensional center manifold. The associated system of ODEs on the center manifold have cubic nonlinearity and is adequately described by 2 coefficients (rather than 8) called the normal form.
 
In Section \ref{sec:s5}, we adapt a recently developed technique called Sparse Identification of Nonlinear Dynamics (SINDy) \cites{Brunton16, Brunton19, Brunton20} which has demonstrated the ability to recover governing equations of complex dynamical systems. The methods presented in SINDy approach the problem of automating the discovery of dynamic equations that describe natural systems through the lens of sparsity-promoting regression techniques such as Least Absolute Shrinkage and Selection Operator (LASSO) \cite{Tibshirani96}. To lend insight into this process, the SINDy algorithm was applied to simulated data from various ROM models to recover their respective sparse equations which is then used to reconstruct the original system's dynamics.
 
 \section{Viscous Moore-Greitzer Equations}
 \label{sec:s2}
 \subsection{Model and Analysis}
 Turbo-jet engine is comprised of 3 parts: axial flow compressor where air gets compressed, the plenum where the air undergoes combustion and rapidly expands, and the turbine where the air is let out. The flow enters from atmospheric pressure at the inlet duct at the left of the figure, proceeds through the compressor block where the static pressure is increased, enters the outlet duct, and then exits to atmospheric pressure through the downstream turbine's throttle. The compressor is made out of an entrance duct, an inlet guide vane (IGV), multiple stages of stator-rotor pairs, and an exit duct towards the plenum. A stator is a rotary system with static blades and a rotor comprises of revolving blades. 
 \begin{figure}[htb]
\centering
\includegraphics[width=3.5in]{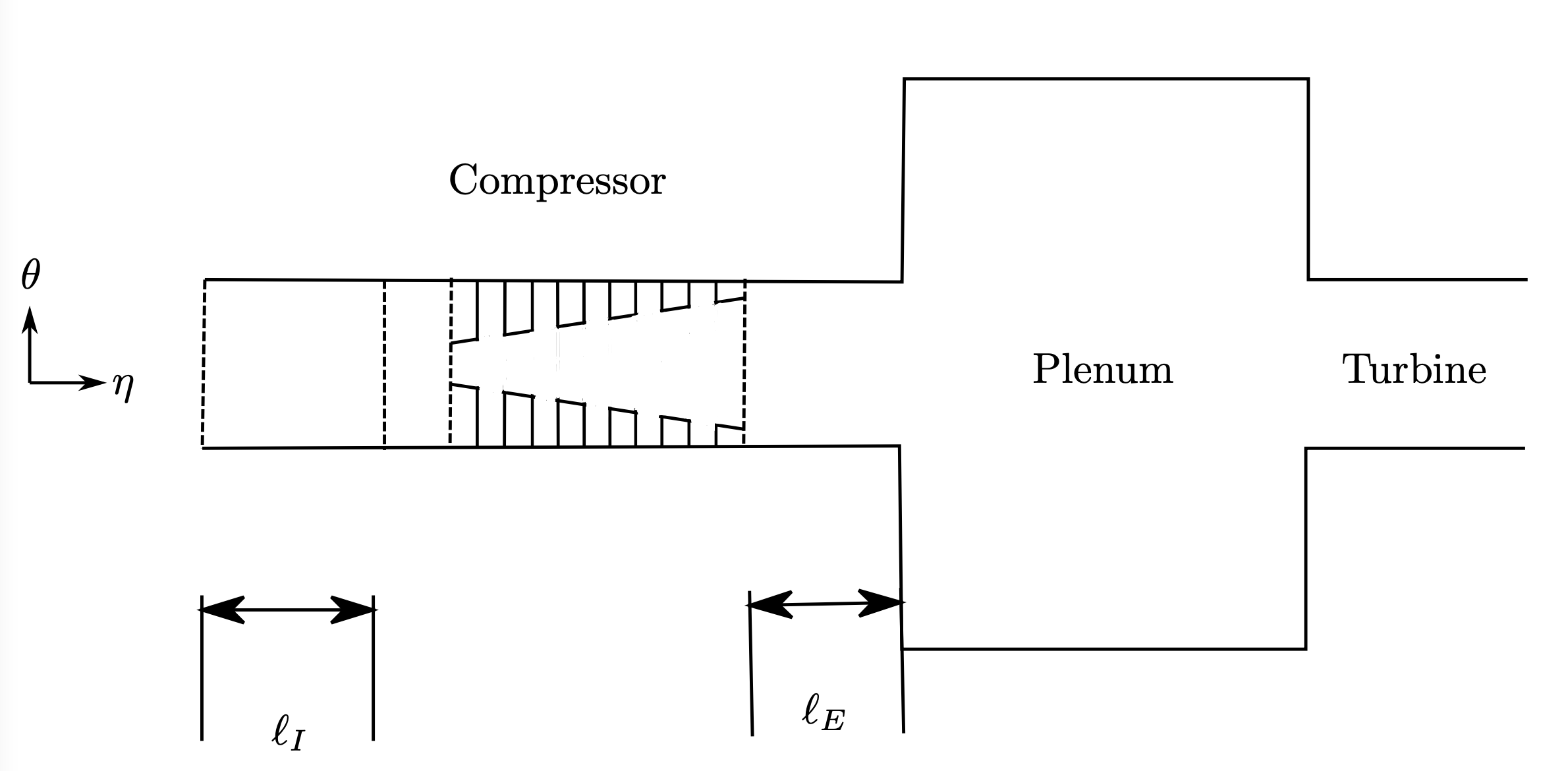}
\caption{Anatomy of a turbo-jet engine comprising of axial flow compressor, plenum, and turbine.}
\label{fig:compressor}
\end{figure}

The basic assumption of the MG compressor model \cites{GM86-1, GM86-2} are:
\begin{enumerate}
\item The pressure rise across the compressor lags behind the pressure drop delivered by the throttle due to mass storage in the exit duct (or plenum).
\item Across the compressor, the difference between the pressure delivered by the compressor and pressure rise that currently exists across the compressor acts to accelerate the flow rate through the compressor.
\item The flow is assumed to be incompressible and irrotational everywhere except inside the plenum where combustion occurs and rapidly expands the air.
\end{enumerate}

The viscous MG equations for a cylindrical axial flow compressor consist of Laplace's partial differential equation (PDE) for disturbance velocity potential $\tilde{\phi}'(t,\theta,\eta)$
\begin{equation}
\tilde{\phi}'_{\eta \eta} + \tilde{\phi}'_{\theta \theta} = 0.
\label{laplace-disturbance}
\end{equation}
with boundary conditions 
\begin{eqnarray}
\psi_c(\Phi(t) + (\tilde{\phi}'_{\eta})_0 ) -  \frac{1}{2 \pi}  \int_0^{2 \pi} \psi_c(\Phi(t) + (\tilde{\phi}'_{\eta})_0 ) d\theta
\nonumber \\
- m (\tilde{\phi}'_{t})_0 - \frac{1}{a} (\tilde{\phi}'_{t \eta})_0 - \frac{1}{2a}(\tilde{\phi}'_{\eta \theta})_0 - \frac{\nu}{2a} (\tilde{\phi}'_{\eta \theta \theta})_0 = 0
\label{bc-eta-0}
\end{eqnarray}
 at $\eta = 0$ and $\tilde{\phi}' = 0$ at $\eta = -\infty$
and a pair of ordinary differential equations (ODEs) for annulus average of axial momentum $\Phi(t)$
\begin{equation}
\Psi(t) + \ell_c \frac{d \Phi(t)}{d t} = \frac{1}{2 \pi}  \int_0^{2 \pi} \psi_c(\Phi(t) + (\tilde{\phi}'_{\eta})_0 ) d\theta.
\label{phi-eqn-viscous-avg}
\end{equation}
and pressure drop from across the compressor $\Psi(t)$
\begin{equation}
\frac{d \Psi(t)}{d t} = \frac{1}{4 B^2 \ell_c} (\Phi(t) - F_T^{-1}(\Psi(t))).
\label{pressure-ode}
\end{equation}
The subscripts of $\tilde{\phi}'$ indicate partial derivatives with respect to time $t$, angular $\theta$ and axial $\eta$ coordinates of the cylindrical compressor. $(\cdot)_0$ means the quantity is evaluated at the compressor entrance $\eta=0$. $a$ is the internal compressor lag, $l_c = l_I + l_E + \frac{1}{a}$ is the characteristic compressor length (dimensionless quantity normalized with respect to compressor radius, see Figure \ref{fig:compressor}), and $B$ is the plenum to compressor volume ratio \cite{Greitzer76}. Detailed derivation of the non-viscous model can be found in \cites{Greitzer76, GM86-1, GM86-2, Moore84-1, Moore84-2, Moore84-3} while the viscous model was developed in \cites{AA93, Mezic98} and thoroughly derived in \cite{BHW07}.

The compressor $\psi_c(\phi)$ and throttle $F_T(\phi)$ characteristic functions that are considered follow \cites{GM86-1, GM86-2}
\begin{eqnarray}
 \psi_c(\phi) &=& \psi_{c_0} + H \left[1 + \frac{3}{2} \left(\frac{\phi}{W} - 1\right) - \frac{1}{2} \left(\frac{\phi}{W} - 1\right)^3 \right]
\label{compressor-characteristic}
\\
 F_T(\phi) &=& \frac{\phi^2}{\gamma^2}.
 \label{throttle-characteristic}
\end{eqnarray}
$H$ and $W$ are the characteristic height and width of the compressor and $\psi_{c_0}$ is a value determined by experiments. Throttle coefficient $\gamma$ describes the amount of opening - large $\gamma$ implies a wide open throttle while small $\gamma$ implies a closed throttle.

Equations \eqref{laplace-disturbance}, \eqref{bc-eta-0}, \eqref{phi-eqn-viscous-avg}, and \eqref{pressure-ode} can be combined into a compact state-space form $\frac{\partial {\bf y}}{\partial t} = {\bf A}{\bf y} + {\bf f}({\bf y})$ following \cite{BH00}
\begin{eqnarray}
\frac{\partial}{\partial t}
\left[
\begin{array}{c}
     g \\ \Phi \\ \Psi
\end{array}
\right] &=&
\left[
\begin{array}{ccc}
     K^{-1} \left(\frac{\nu}{2} \frac{\partial^2}{\partial \theta^2} - \frac{1}{2} \frac{\partial}{\partial \theta}\right) & 0 & 0
     \\
     0 & 0 & 0
     \\
     0 & 0 & 0
\end{array}
\right]
\left[
\begin{array}{c}
     g \\ \Phi \\ \Psi
\end{array}
\right]
\nonumber \\
&&+
\left[
\begin{array}{c}
     a K^{-1} (\psi_c(\Phi + g) - \bar{\psi}_c)
     \\
     \frac{1}{l_c} (\bar{\psi}_c - \Psi)
     \\
     \frac{1}{4B^2 l_c} (\Phi - \gamma \sqrt{\Psi})
\end{array}
\right]
\label{ayplusf}
\end{eqnarray}
by introducing state variable $g$
\begin{equation}
     g(t,\theta) = (\tilde{\phi}'_{\eta})_0 = \sum_{n \in Z} |n| \tilde{\phi}'_n(t) e^{in \theta} = \sum_{n \in Z}  g_n e^{in \theta}.
     \label{g-sum}
\end{equation}
where
\begin{equation}
    \tilde{\phi}'(t, \theta, \eta) = \sum_{n \in Z} \tilde{\bar{\phi}}'_n(t) e^{|n|\eta + in \theta},
\end{equation}
is the solution to \eqref{laplace-disturbance} and we define
\begin{equation}
    \bar{\psi}_c = \frac{1}{2\pi} \int_0^{2\pi} \psi_c(\Phi + g) d\theta.
    \label{psi-c-bar}
\end{equation}
as well as an operator $K$ that acts on $\phi = \sum_{n \in Z} \tilde{\phi}_n e^{in \theta}$ such that 
\begin{equation}
    K(\phi) = \sum_{n \in Z} \left(1 + \frac{ma}{|n|}\right)\tilde{\phi}_n e^{in \theta}.
    \label{K-operator}
\end{equation}

To inspect the nonlinearities in ${\bf f}({\bf y})$, we perform Taylor series' expansion on $\psi_c(\Phi+g)$ up to the third cubic term to expand the integrand of $\bar{\psi}_c$
\begin{eqnarray}
\bar{\psi_c} &=& \psi_c(\Phi) + \frac{1}{2} \psi_c''(\Phi) \sum_{m, n \in Z}^{m+n = 0} g_m g_n
\nonumber \\
&&+ \frac{1}{6} \psi'''_c(\Phi) \sum_{k, m, n \in Z}^{k+m+n = 0} g_k g_m g_n.
\end{eqnarray}
Note that $g(t,\theta)$ has a vanishing average property due to assumptions made to the disturbance flow. Therefore, $\bar{\psi}_c$ is only a function of $t$ and not $\theta$ and as a result, $K^{-1}(\bar{\psi}_c) = 0$. The nonlinearity vector ${\bf f}({\bf y})$ becomes
\begin{eqnarray}
{\bf f}({\bf y}) = \left[
\begin{array}{c}
a K^{-1} (\psi_c'(\Phi) g + \frac{1}{2} \psi''_c(\Phi) g^2 + \frac{1}{6} \psi'''_c(\Phi) g^3)
     \\
     \frac{1}{l_c} (\psi_c(\Phi) + \frac{1}{2} \psi_c''(\Phi) \sum_{m, n \in Z}^{m+n = 0} g_m g_n \hdots
     \\
     + \frac{1}{6} \psi'''_c(\Phi) \sum_{k, m, n \in Z}^{k+m+n = 0} g_k g_m g_n - \Psi)
     \\
     \frac{1}{4B^2 l_c} (\Phi - \gamma \sqrt{\Psi})
\end{array}
\right].
\label{f-expanded}
\end{eqnarray}

The system \eqref{ayplusf}'s equilibrium consist of $g_e(\theta) = 0$ and $\Psi_e = \psi_c(\Phi_e) = F_T(\Phi_e)$ which means $(\Phi_e, \Psi_e)$ lies on the intersection of curves \eqref{compressor-characteristic} and \eqref{throttle-characteristic}. $\Phi_e$ can be solved by finding the root of the polynomial
\begin{equation}
    -\frac{H}{2W^3} \Phi_e^3 + \left(
    \frac{3H}{2W^2} - \frac{1}{\gamma^2}
    \right) \Phi_e^2 + \psi_{c_0} = 0.
    \label{equilibrium-poly}
\end{equation}
\eqref{equilibrium-poly} has one real root and a pair of imaginary roots, where the real root is
\begin{eqnarray}
    \Phi_e &=& \sqrt[3]{X - Y^3 + \sqrt{X(X - 2Y^3)}}
    \nonumber \\
    &&+ \sqrt[3]{X - Y^3 - \sqrt{X(X - 2Y^3)}} - Y
    \label{equilibrium-phi-explicit}
\end{eqnarray}
and
\begin{equation}
    X = \frac{W^3}{H}\psi_{c_0}, \qquad
    Y = \frac{2W^3}{3H} \left(\frac{1}{\gamma^2} -
    \frac{3H}{2W^2}\right).
    \label{XandY}
\end{equation}
For our analysis, $\gamma$ is the bifurcation parameter to be varied for different kinds of Hopf bifurcation. 

The Jacobian of ${\bf f}({\bf y})$ at equilibrium is
\begin{eqnarray}
\nabla{\bf f}_{{\bf y}_{e}} = \left[
\begin{array}{ccc}
a K^{-1}(\psi_c'(\Phi_e)) & 0 & 0 \\
0 & \frac{1}{l_c} \psi_c'(\Phi_e) & -\frac{1}{l_c}  \\
0 & \frac{1}{4 B^2 l_c} & - \frac{1}{4 B^2 l_c} \frac{\gamma^2}{2 \Phi_e}
\end{array}
\right].
\end{eqnarray}
The eigenvalues of $({\bf A} + \nabla {\bf f}_{{\bf y}_e})$ corresponding to the PDE are
\begin{equation}
\lambda_n = \left(\frac{a |n|}{|n| + am}\right) \left(\psi_c'(\Phi_e) - \frac{\nu}{2a} n^2 - \frac{1}{2a} (in) \right)
\label{PDEev}
\end{equation}
and the eigenvalues of $({\bf A} + \nabla {\bf f}_{{\bf y}_e})$ corresponding to the ODEs are
\begin{eqnarray}
\mu_{1,2} &=& \frac{1}{2 l_c} \left[\left(\psi_c'(\Phi_e) - \frac{\gamma}{8B^2\sqrt{\Psi_e}}\right) \right.
\nonumber \\
&&\left. \pm \sqrt{\left(\psi_c'(\Phi_e) + \frac{\gamma}{8B^2\sqrt{\Psi_e}}\right)^2 - \frac{1}{B^2}} \right].
\label{ODEev}
\end{eqnarray}

Hopf bifurcation occurs when a pair of $({\bf A} + \nabla {\bf f}_{{\bf y}_e})$ eigenvalues' real parts cross the imaginary axis with the derivative of the real parts with respect to $\gamma$ is not equal to zero. There are three possibilities: surge (ODE bifurcation), stall (PDE bifurcation), and combination (simultaneous ODE and PDE bifurcations).

The critical bifurcation point for surge is $\gamma_{c,surge}$ such that $Re({\mu}_{1,2}) = 0$. When $\gamma < \gamma_{c,surge}$, surge occurs. It is difficult to obtain an explicit expression for $\gamma_{c,surge}$ but $\gamma_{c,surge}$ is the solution to
\begin{equation}
\Phi_e(\gamma_{c,surge}) \left(2 - \frac{\Phi_e(\gamma_{c,surge})}{W} \right) - \frac{\gamma_{c,surge}^2}{4B^2} \frac{W^2}{3H}= 0.
\end{equation}
The condition for surge is $\left. \frac{\partial}{\partial \gamma} (Re(\mu_{1,2})) \right|_{\gamma_{c,surge}} > 0$.

The critical bifurcation point for stall is $\gamma_{c,stall}$ such that $Re({\lambda}_1) = 0$. When $\gamma < \gamma_{c,stall}$, stall occurs. Again, it is difficult to obtain an explicit expression for $\gamma_{c,stall}$ but $\gamma_{c,stall}$ is the solution to
\begin{equation}
\Phi_e(\gamma_{c,stall}) \left(2 - \frac{\Phi_e(\gamma_{c,stall})}{W} \right) - \frac{\nu W^2}{3aH}= 0.
\end{equation}
The condition for stall is $\left. \frac{\partial}{\partial \gamma} (Re(\lambda_1)) \right|_{\gamma_{c,stall}} > 0$.

It is possible for the largest PDE eigenvalue pairs and both ODE eigenvalues to simultaneously cross the imaginary axis. This is achieved when $\gamma = \gamma_{c,combo}$ where
\begin{equation}
\psi_c'(\Phi_e(\gamma_{c,combo})) = \frac{\gamma_{c,combo}^2}{8B^2} \frac{1}{\Phi_e(\gamma_{c,combo})} = \frac{\nu}{2a}.
\end{equation}
For the combination case, it is possible to calculate the expression for the normal form which are the diagonal entries of $({\bf A} + \nabla {\bf f}_{{\bf y}_e})(\gamma_{c,combo})$
\begin{eqnarray}
{\bf D}(\gamma_{c,combo}) &=& {\bf T}^{-1}({\bf A} + \nabla {\bf f}_{{\bf y}_e}){\bf T}(\gamma_{c,combo})
\nonumber \\
&=& \text{diag} \left(
\left[
\begin{array}{c}
(K^{-1} \left(\frac{\nu}{2}(1-n^2) - \frac{1}{2} ni \right)
\\
\frac{1}{2 l_c} \sqrt{\frac{\nu^2}{a^2} - \frac{1}{B^2}}
\\
-\frac{1}{2 l_c} \sqrt{\frac{\nu^2}{a^2} - \frac{1}{B^2}}
\end{array}
\right]
\right).
\end{eqnarray}

\subsection{Rotating Stall Simulation}
The system of equations \eqref{ayplusf} is integrated using the spectral method. $\theta \in [-\pi, \pi)$ is discretized into 512 equally spaced points, leading to a system of 514 ODEs (512 of which are Fourier coefficients of $g(t,\theta)$) to be numerically integrated using SciPy's \texttt{solve\_ivp} with $dt=0.1$. The following parameter values are used in all cases
\begin{eqnarray}
l_c &=& 8, \quad m = 1.75, \quad a = 1/3.5, \quad \nu = 1
\nonumber \\
\psi_{c_0} &=& 1.67 H, \quad H = 0.18, \quad W = 0.25.
\end{eqnarray}
The plenum to compressor volume ratio $B$ and the throttle opening $\gamma$ are chosen to produce different type of bifurcations.
\begin{figure}[H]
\includegraphics[width=3.5in]{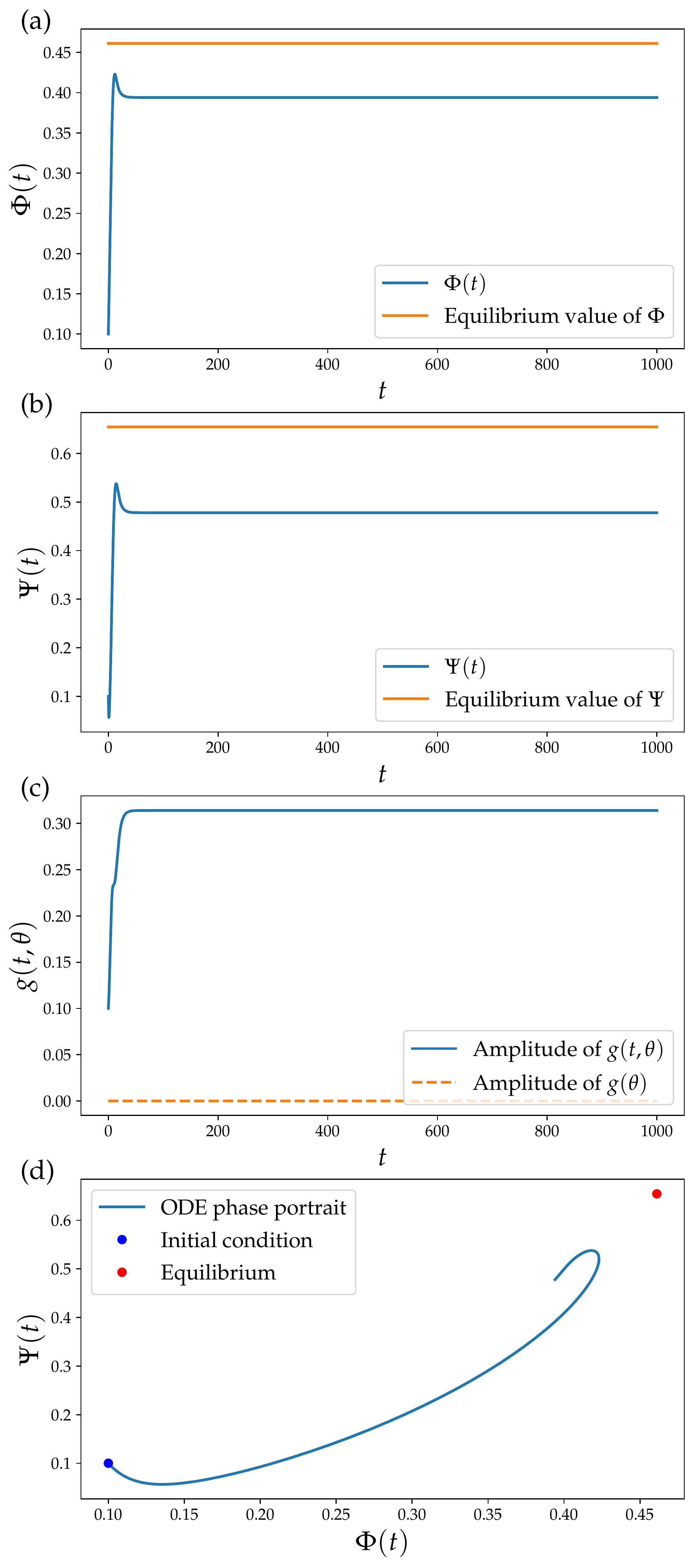}
\caption{Stall dynamics of viscous MG equations with $B = 0.15$ and $\gamma = 0.57$. $\lambda_{1,-1} = 0.077 \mp i$ and $\mu_{1,2} = -0.23 \pm 0.32i$. (a), (b) $\Phi(t)$ and $\Psi(t)$ do not settle at their stable equilibrium values due to influence from PDE Hopf bifurcation. (c) Amplitude of $g(t,\theta)$ in $t$ settles to a non-zero value during Hopf bifurcation. (d) Phase portrait of ODE states $\Phi(t)$ and $\Psi(t)$ which does not settle at the $(\Phi_e, \Psi_e)$, but at $(\Phi_{ss}, \Psi_{ss})$.}
\end{figure}

Of particular interest is the simulation of the stall case which corrects \cite{Xiao08} as we carefully incorporate the quadratic and cubic terms in \eqref{f-expanded} that do not vanish in the simulation. We observe a standing wave limit cycle in the PDE solution which causes the ODE solutions to not settle at their equilibrium values at steady state. If we approximate the long term behavior of the PDE as $g_{ss}(\theta) = A \cos({\theta})$, plug in this assumption to \eqref{phi-eqn-viscous-avg}, and set $\frac{d\Phi}{dt} = \frac{d\Psi}{dt} = 0$, $\Phi_{ss}$ can be found by solving for the root of
\begin{eqnarray}
\psi_c(\Phi_{ss}) + \frac{1}{2} \psi_c''(\Phi_{ss}) \frac{A^2}{2} - \frac{\Phi^2_{ss}}{\gamma^2} &=& 0.
\end{eqnarray}
Subsequently $\Psi_{ss} = \frac{\Phi_{ss}^2}{\gamma^2}$.

\section{Reduced Order Modelling (ROM)}
\label{sec:s3}
\subsection{Principal Component Analysis (PCA)}
PCA is a method to find principal axes in high dimensional data. These principal axes span the eigenvectors of the covariance matrix of the measurements which are orthonormal to each other such that the individual data along these directions are linearly uncorrelated. PCA can also be used as a dimensional reduction tool by truncating a measurement's linear combination in its principal axes. Constructing basis functions from data using PCA can be formulated mathematically as a low-rank matrix approximation problem which can be easily computed by using the singular value decomposition (SVD) \cite{GK64}. PCA is also known as different names such as proper orthogonal decomposition (POD) in mechanical engineering \cite{HLBR12} and discrete Karhunen-Lo{\`e}ve expansion in signal processing and information theory \cites{Karhunen47, Loeve78}. Our work is inspired by Karhunen-Lo{\`e}ve expansion applied to find reduced dynamics of turbulent flows \cite{Sirovich87} and atmospheric waves \cite{CZN07}.

Suppose we have $N$ observations of $n$-dimensional data
\begin{equation}
{\bf Y} = [{\bf y}_1 \hdots {\bf y}_N]
\end{equation}
where ${\bf Y} \in \mathbbm{R}^{n \times N}, {\bf y_i} = {\bf y}(t_i) \in \mathbbm{R}^{n}$. After centering the data about its empirical mean to get ${\bf Y}_0$, define a transformation ${\bf W} \in \mathbbm{R}^{n\times m}$ where $m< n$ in the context of dimensional reduction. The lower dimensional data is calculated by
\begin{equation}
    {\bf X}_0 = {\bf W}^T {\bf Y}_0.
\end{equation}
If we use PCA, then the transformation is defined as
\begin{equation}
    {\bf W} = {\bf P}_m = [{\bf p}_1 \hdots {\bf p}_m]
    \label{pcatransform}
\end{equation}
where ${\bf p}_m$ are the principal axes of ${\bf Y}_0$ or the first $m$ eigenvectors of the covariance matrix ${\bf Y}_0 {\bf Y}_0^T$.

\subsection{Neural Network Implementation of PCA}
The most widely known neural network architecture is the multilayer feedforward neural network (FNN) which is also known as multilayer perceptron (MLP). A multilayer FNN consists of a number layers starting with an input layer followed by one or more hidden layers and ending with an output layer all are connected in feedforward manner. 

\begin{figure}[htb]
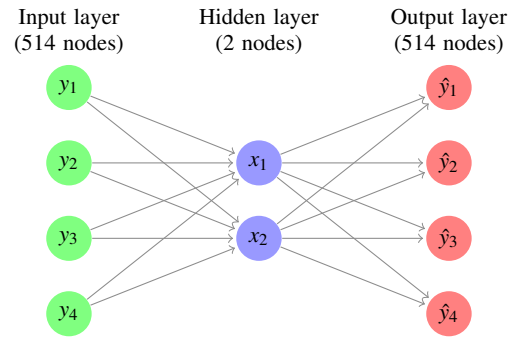

\begin{center}
\begin{neuralnetwork}[height=4]
        \newcommand{\x}[2]{$y_#2$}
        \newcommand{\y}[2]{$\hat{y}_#2$}
        \newcommand{\zfirst}[2]{\small $x_#2$}
        \inputlayer[count=4, bias=false, title=Input layer \\ (514 nodes), text=\x]
        \hiddenlayer[count=2, bias=false, title=Hidden layer \\ (2 nodes), text=\zfirst]
        \linklayers
        \outputlayer[count=4, title=Output layer \\ (514 nodes), text=\y] \linklayers
    \end{neuralnetwork}
\caption{The simplest neural network architecture}    \label{fig:Autoencoder1}  
\end{center}
\end{figure}

An autoencoder is a type of multilayer feedforward neural network that at its simplest form (as illustrated in Figure \ref{fig:Autoencoder1})has an input layer with $n$ nodes, followed by a hidden layer with $m$ nodes (where $m<n$), followed by an output layer with $n$ nodes. When the activation functions are chosen to be linear, the input-output relationship is given by
\begin{equation}
    {\bf \hat{y}} = {\bf W}_2 ({\bf W}_1 {\bf y} + {\bf b}_1) + {\bf b}_2
\end{equation}
where ${\bf W}_1, {\bf W}_2^T \in \mathbbm{R}^{m \times n}$ are the encoder and decoder weight matrices, and ${\bf b}_1 \in \mathbbm{R}^{m}$, ${\bf b}_2 \in \mathbbm{R}^{n}$ are the encoder and decoder bias vectors.
Once the optimal $\{{\bf W}_1, {\bf W}_2, {\bf b}_1, {\bf b}_2\}$ are found, we can construct an encoder to reduce the input into a reduced order data ${\bf z} \in \mathbbm{R}^m$ using $\{{\bf W}_1, {\bf b}_1\}$ and a decoder to convert the encoded data back to its original dimension using $\{{\bf W}_2, {\bf b}_2\}$.

Under certain assumptions on the error function landscape, the minimization problem for the autoencoder reduces to
\begin{equation}
    \min_{\{{\bf W}_2\}} \lVert {\bf Y}_0 - {\bf W}_2{\bf W}_2^+ {\bf Y}_0 \rVert_F^2.
    \label{pca-nn-loss}
\end{equation}
where ${\bf W}_2^+$ is the Moore-Penrose inverse/pseudoinverse \cites{Moore20, Penrose55} of ${\bf W}_2$.
For the case when the columns of ${\bf W}_2$ are orthonormal like ${\bf P}_m$, then ${\bf W}_2^+ = {\bf W}_2^T$ will make \eqref{pca-nn-loss} equal to the reconstruction error of PCA. Therefore, it is clear that ${\bf P}_m$ is a solution to the autoencoder optimization problem \cite{EY36}. The problem is that the product of ${\bf P}_m$ with any proper orthogonal matrix ${\bf Q} \in \mathbbm{R}^{m\times m}$ will be a minimizer ${\bf W}_2$, such that there are infinitely many solutions. Coupled with the fact that mini-batch stochastic gradient descent \cite{LZCS14} is the go-to optimization algorithm in today's neural network frameworks, there is no guarantee that ${\bf W}_2$ converges to the same value when the training procedure is repeated, let alone align itself to ${\bf P}_m$. While any ${\bf W}_2$ in this space can be used to mimic the input data almost perfectly, this inconsistency is an issue in our problem as we would like to further uncover the underlying structure of the encoded measurements ${\bf X} = [{\bf x}_1 \hdots {\bf x}_N]$.

\subsubsection{Regularized Linear Autoencoder}
An approach to recover the PCA principal axes from autoencoder weights is based on the following hypothesis \cite{Plaut18}: the first $m$ left singular vectors of ${\bf W}$ is also the first $m$ principal axes of ${\bf Y}_0$. This hypothesis can be framed as an autoencoder with a regularizer or penalty to the sum of the Frobenius norms of the encoder weight matrix ${\bf W}_1$ and decoder weight matrix ${\bf W}_2$ 
\begin{equation}
    \min_{\{{\bf W}_1, {\bf W}_2\}} \lVert {\bf Y}_0- {\bf W}_2{\bf W}_1 {\bf Y}_0 \rVert_F^2 + \lambda (\lVert {\bf W}_1\rVert_F^2 + \lVert{\bf W}_2 \rVert_F^2).
\end{equation}
For a large enough $\lambda$ value, the error surface is guaranteed to be convex with a single global minima which will correspond to the PCA principal axes \cite{KBGS19}. Additionally, the minimum values of this loss function is ${\bf W}_1^* = {\bf W}_2^T$ unlike ${\bf W}_1^* = {\bf W}_2^+$ in the original approach. ${\bf W}_2$ is also found to be equal to the principal axes of probabilistic PCA \cite{TB99} when $\sigma^2 = \lambda$, $\sigma$ being the variance of the data in the Bayesian framework/maximum aposteriori estimation (MAP) derivation of probabilistic PCA.

\subsubsection{Nonlinear Principal Component Analysis (NLPCA) and Autoencoder}
\begin{figure*}
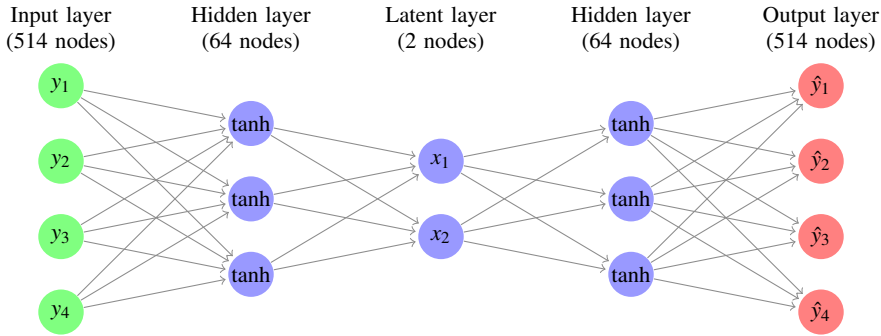

\begin{center}
\begin{neuralnetwork}[height=4]
        \newcommand{\x}[2]{$y_#2$}
        \newcommand{\y}[2]{$\hat{y}_#2$}
        \newcommand{\hfirst}[2]{\small $\tanh$}
        \newcommand{\zfirst}[2]{\small $x_#2$}
        \newcommand{\hsecond}[2]{\small $\tanh$}
        \inputlayer[count=4, bias=false, title=Input layer \\ (514 nodes), text=\x]
        \hiddenlayer[count=3, bias=false, title=Hidden layer \\ (64 nodes), text=\hfirst] \linklayers
        \hiddenlayer[count=2, bias=false, title=Latent layer \\ (2 nodes), text=\zfirst]
        \linklayers
        \hiddenlayer[count=3, bias=false, title=Hidden layer \\ (64 nodes), text=\hsecond] \linklayers
        \outputlayer[count=4, title=Output layer \\ (514 nodes), text=\y] \linklayers
    \end{neuralnetwork}
\caption{Autoencoder architecture for MG compressor data with $tanh()$ activations function}    \label{fig:Autoencoder2}  
\end{center}
\end{figure*}

NLPCA was developed to uncover the underlying nonlinear manifold in large dimensional datasets. It was first implemented using neural network in \cite{Kramer91}. The neural network architecture we are considering to train our NLPCA autoencoder is shown in Figure \ref{fig:Autoencoder2}
\begin{eqnarray}
{\bf X} &=& {\bf W}_2 \tanh({\bf W}_1 {\bf Y} + {\bf b}_1) + {\bf b}_2
\nonumber \\
\hat {\bf Y} &=& {\bf W}_4 \tanh({\bf W}_3 {\bf X} + {\bf b}_3) + {\bf b}_4.
\label{nlpca-nn}
\end{eqnarray}
We choose the nonlinear activation function $\tanh()$ as in \cite{Scholz02} under the justification that a trigonometric function would fit well with the solutions of the MG equations which are spanned by the Fourier basis \eqref{g-sum}. The NLPCA autoencoder is trained to minimize the loss function of
\begin{equation}
    \min_{\{{\bf W}_{1,2,3,4}, {\bf b}_{1,2,3,4}\}} \lVert {\bf Y} - \hat{\bf Y}\rVert_F^2.
\end{equation}
The resulting $\{{\bf W}_{1,2,3,4}\}$ and $\{{\bf b}_{1,2,3,4}\}$ are then used to construct an encoder and decoder as per \eqref{nlpca-nn}.

\section{Sparsity in Reduced Order Data}
\label{sec:s4}
Over the past two decades, researchers have focused on sparsity as one type of low-dimensional structure. Given the recent advances in both compressed sensing \cites{CRT06, Donoho06, CW08} and sparse regression \cite{Tibshirani96}, it has become computationally feasible to extract system dynamics from large multimodal datasets. The term sparse in signal processing context refers to the case where signals (or any type of data, in general) have few non-zero components with respect to the total number of components. 
It is well known in dynamical systems, the normal forms provide a way of finding a coordinate system in which the dynamical system takes the ``simplest" or  ``minimal" form. 
The normal forms, which are sparse in the space of homogeneous vector polynomial of certain degree,  is  calculated by making judicious choices of the solutions to the homological equations~\cite{GH83}. 
Hence, in the context of our work, close to the bifurcation point, the sparse regression techniques rely heavily on the fact that many dynamical systems can be represented by governing equations that are sparse in the space of all possible functions of a given algebraic structure.

\subsection{Compressed Sensing}
Compressed sensing (CS) is a technique for sampling and reconstructing sparse signals, i.e. signals that can be represented by $k<<n$ significant coefficients over an $n$- dimensional basis. The central goal of CS is the recovery of sparse vectors from a small number of linear measurements, which distinguishes CS from other dimensionality reduction techniques. Hence, this allows for polynomial-time reconstruction of the sparse signal \cite{Donoho06}.

In \cite{Donoho06} and \cite{CRT06}, the original sparse ($k$-sparse) signal is projected onto a lower-dimensional subspace via a random projection scheme, called the sampling matrix. More precisely, this broader objective is exemplified by the important special case in which one is interested in finding a vector $X \in \mathbbm{R}^n$ using the (noisy) observation or the measurement data 
\begin{equation}
Y = {\Theta} X + {\eta}, \quad\text{where } {\Theta}\in\mathbbm{C}^{m\times n}\;\;\text{with}\;\; k<m<n,
\end{equation}
is the known sensing or sampling matrix and $\eta$ is the measurement noise.

In general, the problem cannot be solved uniquely. However, if $X$ is $k$-sparse i.e., if it has up to $k$ non-zero entries, the theory of CS shows that it is possible to reconstruct $X$, a $k$-sparse vector in $\mathbbm{R}^n$ uniquely from $m$ linear measurements even when $m<<n$, by exploiting the sparsity of $X$. This can be achieved by finding the sparsest signal consistent with the vector of measurements \cite{Donoho06}, i.e.
\begin{equation}
\argmin_{X\in \mathbbm{R}^{n}} \;
\lVert X \rVert_{0} \;\; \text{subject to } \lVert Y -\Theta X \rVert_{_{2}}\leq \varepsilon
\label{P:l0}
\end{equation}
where $\lVert X \rVert_{0}$ denotes the $l_{0}$ norm for $X$ (the number of non-zero entries of $X$), while $\varepsilon$ denotes a parameter that depends on the level of measurement noise $\eta$.
It can be shown that the $l_{0}$  minimization method can exactly reconstruct the original signal in the absence of noise using a properly chosen sensing matrix $\Theta$ whenever $m>2k$. However, $l_{0}$ minimization problem \eqref{P:l0} is a non-convex problem which is NP-hard.

Instead of problem \eqref{P:l0} we consider its $l_{1}$ convex relaxation which may be stated as \cite{CDS98}
\begin{equation}
\argmin_{X \in \mathbbm{R}^{n}} \;
\lVert X \rVert_{1} \;\; \text{subject to } \lVert Y -\Theta X \rVert_{_{2}}\leq \varepsilon
\label{P:l1}
\end{equation}
where the $l_1$ norm (sum of the absolute values of the entries of $X$) is a convex function. Hence \eqref{P:l1} is a convex optimization problem which can accurately approximate the solution to \eqref{P:l0} in polynomial time with high probability if measurement matrix $\Theta$ is chosen to satisfy a necessary condition called ``Restricted Isometry Property" (RIP) \cites{CT05, CRT06}. Loosely speaking, if $\Theta$ satisfies the RIP condition, then the measurement matrix approximately preserves the Euclidean length of every $k$-sparse signal. Equivalently, all subsets of $k$ columns taken from $\Theta$  are nearly orthogonal. One should note that the $l_{1}$ minimization in \eqref{P:l1} is closely related to the LASSO problem \cite{Tibshirani96}
\begin{equation}
  \argmin_{X \in \mathbbm{R}^{n}}
  \lVert Y - \Theta X\rVert_2^2 + \alpha \lVert X \rVert_{1}
  \label{P:l2}
\end{equation}
where $\alpha \geq 0$ is a regularization parameter. If $\varepsilon$ and $\alpha$ in \eqref{P:l1} and \eqref{P:l2} satisfy some special conditions, the two problems are equivalent; however, characterizing the relationships between $\varepsilon$ and $\alpha$ is difficult except for the special case of orthogonal sensing matrices $\Theta$. The practical success and importance of the lasso can be attributed to the fact that in many cases $X$ is sparse.

\subsection{Sparse Identification of Dynamical Systems (SINDy)}
Sparse identification of nonlinear dynamics (SINDy) \cite{Brunton16} is an algorithm for discovering the dynamical equations directly from the data. The problem of model discovery from data can be formulated as a feature selection problem in machine learning \cite{KR92}. The SINDy algorithm takes $m$-time measurements of ${\bf x} \in \mathbbm{R}^n$, ${\bf X} = [{\bf x}(t_1), \hdots, {\bf x}(t_m)]^T \in \mathbbm{R}^{m\times n}$ and attempts to discover the structure of a nonlinear differential equation of the form
\begin{equation}
\dot{{\bf X}} = {\bf f}({\bf X}(t)) \approx {\bTheta}({\bf X}) {\bXi}
\end{equation}
where ${\bTheta}({\bf X})=[\theta_{1}({\bf X}), \theta_{2}({\bf X}),...,\theta_{p}({\bf X})]\in\mathbbm{R}^{m \times p}$ form the dictionary of basis functions, and ${\bXi} \in\mathbbm{R}^{p \times n}$ is the matrix of coefficients, where each column corresponds to an equation with $p$ terms. $p$ is the maximal number of $n$-multivariate monomials of degree at most $d$. The majority of $\bXi$ entries are zero while the remaining non-zero entries identify the active terms contributing to the sparse representation of the dynamics ${\bf f}({\bf X})$. To guarantee sparsity, SINDy is reformulated as a LASSO problem 
\begin{equation}
\argmin_{\bXi} \frac{1}{m} \sum_{i=1}^m \lVert \dot{\bf x}(t_i)-{\bTheta}({\bf x}(t_i)) {\bXi} \rVert_2^2 + \alpha \lVert {\bxi} \rVert_1.
\label{eq:sindy-cost}
\end{equation}
where $\bxi(\bXi) \in \mathbbm{R}^{1 \times pn}$ is a vector of all entries inside $\bXi$. LASSO is an optimization algorithm that finds a sparse solution for \eqref{eq:sindy-cost} by initializing $\bXi = {\bf 0}$ and at each iteration, it tries to find an update for $\bXi$ one matrix entry at a time. The $l_1$ regularization coefficient $\alpha$ acts as a threshold such that if the an optimal condition involving $\alpha$ is not satisfied for a particular $\bXi$ entry, the entry is chosen to be equal to zero. Increasing the value of $\alpha$ leads to more zero entries in ${\bXi}$, resulting in a sparse model.

The dictionary of basis functions for monomial sampling of dynamical system is
\begin{eqnarray}
{\bTheta}({\bf X}) = \left[
\begin{array}{ccccc}
| & | & | & | & | \\
\bf{1} & \bf{X} & \bf{X}^{P_2} & \bf{X}^{P_3} & \hdots\\
| & | & | & | & |
\end{array}
\right].
\label{eq:cubictheta}
\end{eqnarray}
The dictionary ${\bTheta}({\bf X})$ is constructed by appending candidate nonlinear functions of ${\bf X}$ column-wise. Here, higher order polynomials are denoted as ${\bf X}^{P_d}$ where $d$ is the order of the polynomial considered. For example, element 1 is a column-vector of ones, element ${\bf X}$ is as defined above, element ${\bf X}^{P_2}$ is the matrix containing the set of all quadratic polynomial functions of the state vector ${\bf x}$, and is constructed as follows:
\begin{eqnarray}
\bf{X}^{P_2} = \left[
\begin{array}{cccccc}
x_1^{2}(t_1) & x_1(t_1)x_2(t_1) &\hdots & x_2^{2}(t_1) & \hdots & x_n^{2}(t_1) \\
x_1^{2}(t_2) & x_1(t_2)x_2(t_2) &\hdots & x_2^{2}(t_2) & \hdots & x_n^{2}(t_2) \\
\vdots & \vdots & \ddots & \vdots & \ddots & \vdots \\
x_1^{2}(t_m) & x_1(t_m)x_2(t_m) & \hdots & x_2^{m}(t_m) & \hdots & x_n^{2}(t_m)
\end{array}
\right].
\end{eqnarray}
We interpolate the reduced MG simulation data as a dynamical systems with cubic nonlinearity which is up to $\bf{X}^{P_3}$.

\section{Compressed Compressor Analysis}
\label{sec:s5}
We run 10 simulations of the viscous MG equations' stall case for $t \in [0,500]$ with $dt = 0.1$. The initial conditions for $g(t,\theta)$'s amplitude, $\Phi$, and $\Psi$ are drawn from the normal distribution with mean 0.1 and standard deviation 0.05. The first 2000 data points (up to $t=200$) containing the transient dynamics are discarded. This gives us 10 ${\bf Y} \in \mathbbm{R}^{3000 \times 514}$ datasets. We perform $k$-fold cross validation \cite{HTF09} on PCA, regularized autoencoder, and NLPCA autoencoder to find the best ROM parameters to bring down the data dimension to 2. Both autoencoders' training were performed using Adam optimizer \cite{KB14} with learning rate of $10^{-4}$ for 10 epochs of 4 mini-batch size for the regularized linear autoencoder and 20 epochs of 4 mini-batch size for the NLPCA autoencoder.

We encode the 10 datasets using the 3 different encoders to obtain 3 versions of 10 ${\bf X} \in \mathbbm{R}^{3000 \times 2}$. For each group of reduced order/encoded data, we perform a cubic nonlinearity dynamical system identification using PySINDy \cite{Brunton20} paired with LASSO optimizer from Python's \texttt{sklearn} package. We train the 3 groups of 10 ${\bf X}$ datasets in order to find the largest $\alpha$ value which maximizes the accuracy ($R^2$ score of the SINDy regression) using grid search \cite{LCBB07}. Another set of equations that are discovered by larger $\alpha$ values to maximize sparsity which only end up capturing the cubic nonlinearities are presented in Appendix \ref{sec:appendix}, as sparsity is a trade-off of accuracy. After finding the most suitable $\alpha$ for each group, we perform another $k$-fold cross validation to decide on a model that best represent the 10 datasets of each ROM.

The discovered reduced governing equations satisfy the normal form if it is sufficiently described by 4 coefficients $\mu, \omega, b_1,b_2$ up to an acceptable numerical tolerance
\begin{eqnarray}
\label{normalform}
\dot{x_1} &=& \mu x_1 - \omega x_2 + b_1 (x_1^2 + x_2^2) x_1 - b_2(x_1^2 + x_2^2) x_2
\nonumber \\
\dot{x_2} &=& \omega x_1 + \mu x_2 + b_2 (x_1^2 + x_2^2) x_1 + b_1(x_1^2 + x_2^2) x_2.
\end{eqnarray}
When the linear operator is semi-simple~(as in Hopf bifurcations), the correct identification of a normal form depends critically on the null space of the homological operator \cite{GH83}.  The consequence of this fact is quite profound. The nonlinear terms in normal form \eqref{normalform} commutes with the linear term. As a consequence, when the equation is normalized to any finite degree $k=3$ and truncated, it will have symmetries that were not present in the original system.

For the reconstruction, the obtained SINDy equations are integrated using the forward Euler method with a fixed integration time step $dt = 0.1$ to be consistent with the chosen smoothed forward difference differentiation scheme. The global truncation error is then subtracted from the raw numerical integration result to correct the estimate. Lastly, the integrated SINDy data are fed into the decoder of the respective reduction methods to reconstruct the high dimensional time series and compared with the original dataset. The datasets and code used to produce the results in this paper can be accessed at \url{https://github.com/alytjong/compressed-compressor}.

\subsection{PCA and SINDy}
\begin{figure}[htb]
    \centering
    \includegraphics[width=3.5in]{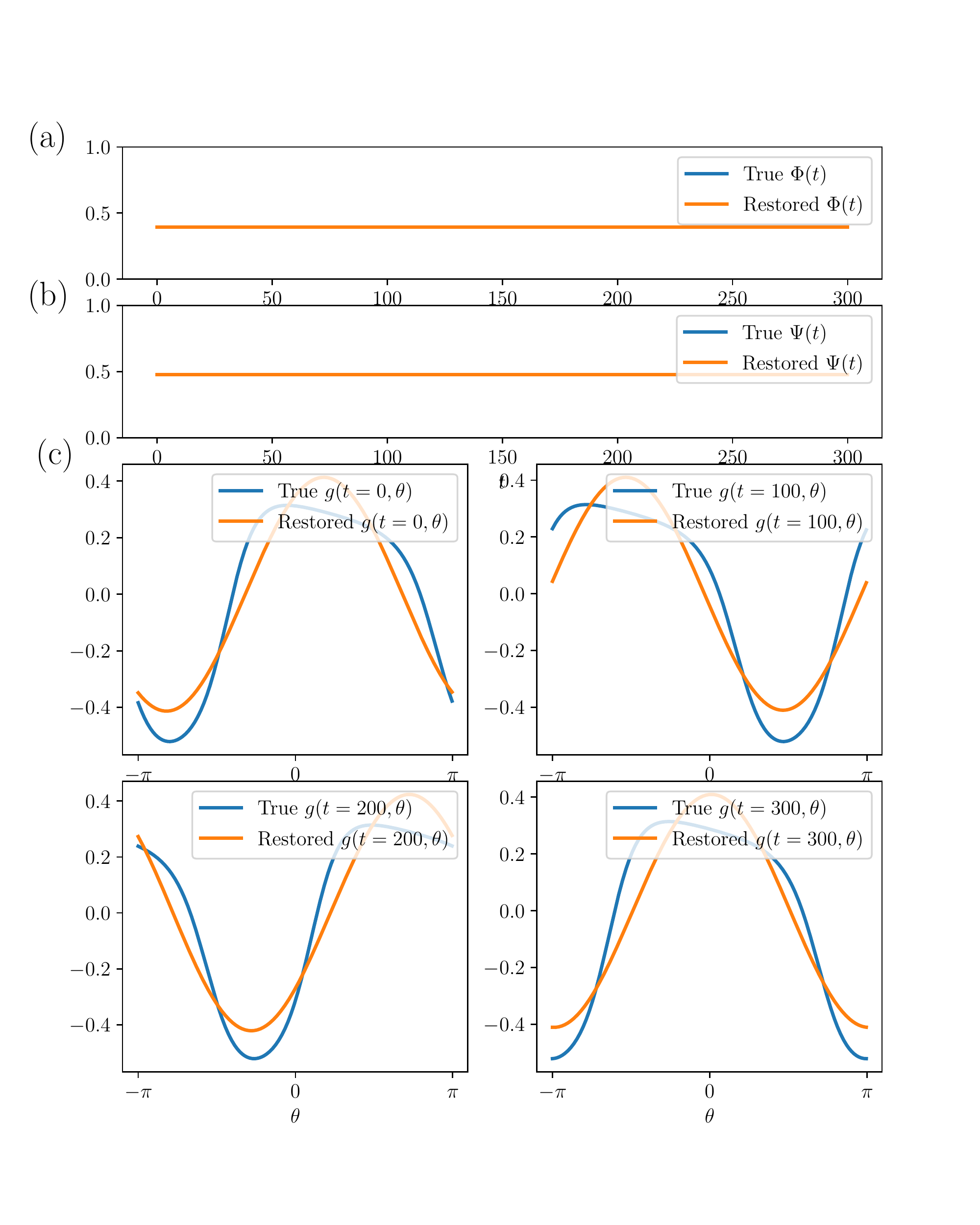}
    \caption{Reconstruction of viscous MG stall dynamics from multiplying ${\bf P}_{m=2}$ with numerical integration result of \eqref{pod-sindy}. (a) $\Phi(t)$ reconstruction result. (b) $\Psi(t)$ reconstruction result. (c) $g(t,\theta)$ reconstruction results at $t=0, 100, 200, 300$.}
    \label{pod-reconstruction}
\end{figure}

The following SINDy regression is obtained using a LASSO threshold of  $\alpha = 0.0035$, which outputs a system of ODEs with 13 coefficients and test $\||{\bf \dot{X}}\||$ score of 0.9999
\begin{eqnarray}
\dot{x}_1 &=& -0.359740 x_2 + 0.000062 x_1^2 -0.000061 x_2^2 
\nonumber \\
&&-0.000144 x_1^3 + 0.000133 x_1^2 x_2 -0.000144 x_1 x_2^2
\nonumber \\
&& + 0.000133 x_2^3
\nonumber \\
\dot{x}_2 &=& 0.312416 x_1 -0.000277 x_1 x_2
\nonumber \\
&&+ 0.000955 x_1^3 -0.000144 x_1^2 x2 + 0.000954 x_1 x_2^2 
\nonumber \\
&&-0.000144 x_2^3.
\label{pod-sindy}
\end{eqnarray}
Some resemblance to the normal form are observed through the almost identical linear frequencies and the repeated cubic coefficient $-0.000144$. Reconstruction result for a chosen random dataset is shown in Figure \ref{pod-reconstruction}.

\subsection{Regularized Linear Autoencoder and SINDy}
\begin{figure}[htb]
    \centering
    \includegraphics[width=3.5in]{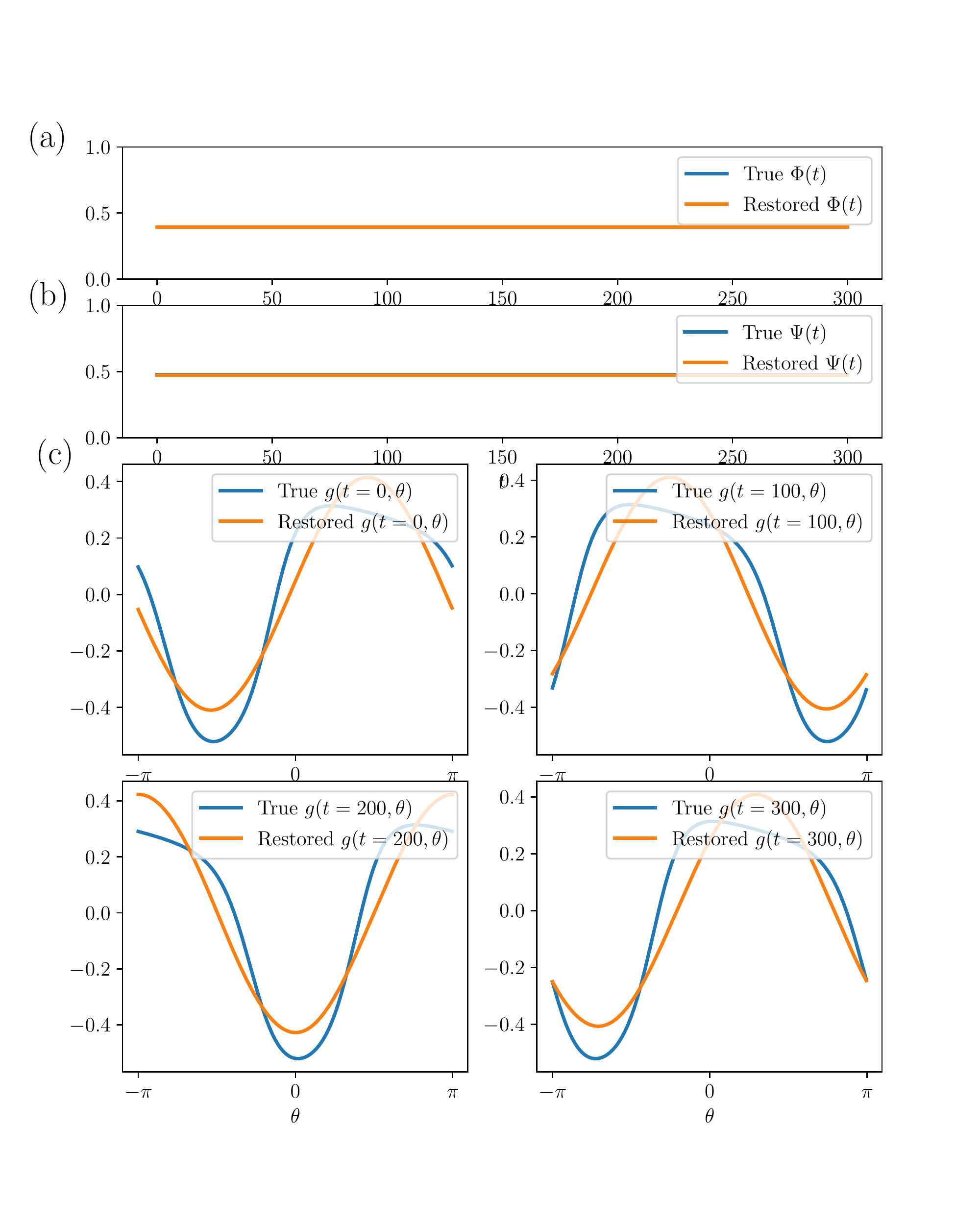}
    \caption{Reconstruction of viscous MG stall dynamics from feeding numerical integration result of \eqref{svdnn-sindy} into the decoder part of a regularized linear autoencoder. (a) $\Phi(t)$ reconstruction result. (b) $\Psi(t)$ reconstruction result. (c) $g(t,\theta)$ reconstruction results at $t=0, 100, 200, 300$.}
    \label{svdnn-reconstruction}
\end{figure}

The following SINDy regression is obtained using a LASSO threshold of  $\alpha = 0.0035$. The output is a system of ODEs with 16 coefficients and test $\||{\bf \dot{X}}\||$ score of 0.9999. A representative equation (since the outcome is always random) is 
\begin{eqnarray}
\dot{x}_1 &=& -0.208500 x_2 + 0.000024 x_1^2 + 0.000099 x_1 x_2 
\nonumber \\
&&-0.000283 x_2^2 -0.000144 x_1^3 -0.003345 x_1^2 x_2
\nonumber \\
&&-0.000145 x_1 x_2^2 -0.003344 x_2^3
\nonumber \\
\dot{x}_2 &=& 0.454352 x_1 + 0.000160 x_1^2 -0.000544 x_1 x_2 
\nonumber \\
&& -0.000119 x_2^2 -0.002309 x_1^3 -0.000144 x_1^2 x_2 
\nonumber \\
&&-0.002308 x_1 x_2^2 -0.000144 x_2^3.
\label{svdnn-sindy}
\end{eqnarray}
Despite the randomness due to different ROM selected at each training, the cubic coefficient $-0.000144$ is always featured. Additionally, while the normal form $\eqref{normalform}$ is not recovered perfectly, some symmetry is still observed in the cubic terms of the individual equations. Reconstruction result for a chosen random dataset is shown in Figure \ref{svdnn-reconstruction}.

\subsection{NLPCA Autoencoder and SINDy}
\begin{figure}[htb]
    \centering
    \includegraphics[width=3.5in]{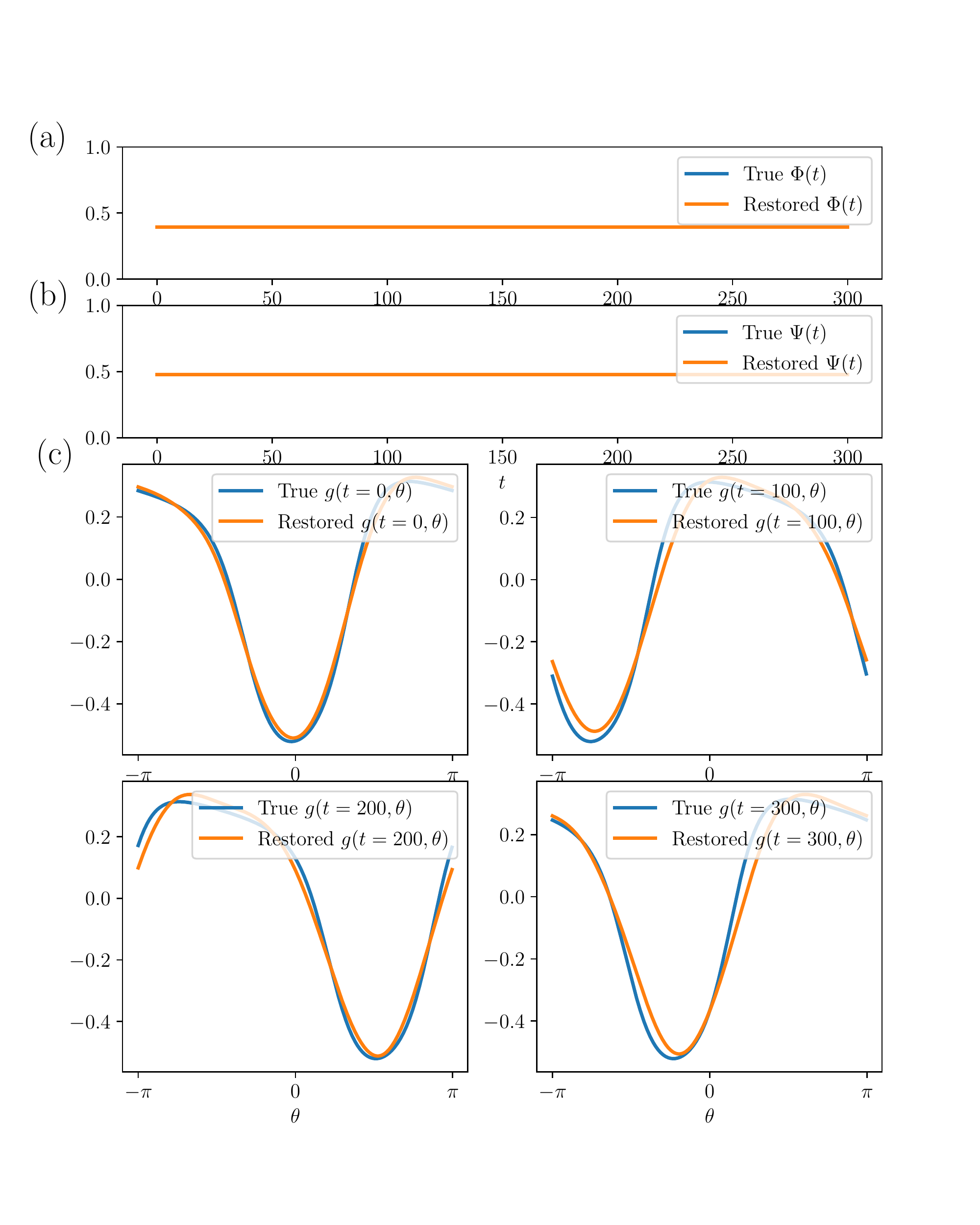}
    \caption{Reconstruction of viscous MG stall dynamics from feeding numerical integration result of \eqref{nlpca-sindy} into the decoder part of an NLPCA autoencoder. (a) $\Phi(t)$ reconstruction result. (b) $\Psi(t)$ reconstruction result. (c) $g(t,\theta)$ reconstruction results at $t=0, 100, 200, 300$.}
    \label{nlpca-reconstruction}
\end{figure}

The following SINDy regression is obtained using a LASSO threshold of  $\alpha = 0.001$. The output is a system of ODEs with 16 coefficients and test $\||{\bf \dot{X}}\||$ score of 0.9999. A representative equation (since the outcome is always random) is 
\begin{eqnarray}
\dot{x}_1 &=& -0.738101 x_2 -0.000244 x_1^2 -0.002102 x_1 x_2 
\nonumber \\
&&+ 0.000428 x_2^2 -0.000038 x_1^3 + 0.007882 x_1^2 x_2
\nonumber \\
&&-0.000319 x_1 x_2^2 + 0.007960 x_2^3
\nonumber \\
\dot{x}_2 &=& 0.603523 x_1 + 0.000096 x_1^2 + 0.000489 x_1 x_2 
\nonumber \\
&&-0.001452 x_2^2 -0.005271 x_1^3 -0.000409 x_1^2 x_2
\nonumber \\
&&-0.004956 x_1 x_2^2 -0.000071 x_2^3.
\label{nlpca-sindy}
\end{eqnarray}
There is no symmetrical structure detected in both equations. Reconstruction result for a chosen random dataset is shown in Figure \ref{nlpca-reconstruction}.

\section{Conclusion and Future Works}
We have showed that it is possible to fully reconstruct the solutions of the viscous MG equations from a system of 2 ODEs up to cubic nonlinearity. It turns out that reconstruction quality is entirely independent of whether the normal form structure of the underlying PDE is detected or not. The NLPCA autoencoder has to be trained for twice as long (double the epoch) compared to the linear autoencoder in order to converge to the local minimum that produces great reconstruction result. Table \ref{table:table-result} summarizes our findings for the three chosen methods.

In order to detect a consistent reduced set of equations representing the PDE Hopf bifurcation, we need the reduced order data to fall along the first two principle axes. Regularization term to linear autoencoder's cost function introduced in \cite{KBGS19} can provide some structure to the discovered SINDy equations, although most of the resulting coefficients will still be random. It would be interesting to find out what regularization term is needed to ensure the convexity of loss landscape of the NLPCA autoencoder in order to obtain both consistent nonlinear structure in the discovered SINDy equations and the best reconstruction result.

\begin{table}
\begin{center}
\begin{tabular}{ | m{1in} | m{0.5in} | m{0.5in} | m{0.5in} | } 
\hline
&PCA & Regulrzd. Linear Autoencoder & NLPCA Autoencoder\\
\hline\hline
Training time  & 25 s & 103 s & 283 s \\
\hline
PDE reconstruction $R^2$ score from training data & 0.8973 & 0.8973 & 0.9916 \\
 \hline
PDE reconstruction $R^2$ score from SINDy equations & 0.8950 & 0.8948 & 0.9887\\
\hline
Number of RHS terms in reduced equations & 13 & 16 & 16\\[1ex] 
\hline
\end{tabular}
\caption{Summary of viscous MG equations' reconstruction from SINDy models identified from different ROMs.}
\label{table:table-result}
\end{center}
\end{table}

Our simple approach rooted in physics-based machine learning which involves a priori knowledge of sparsity and the center manifold theory \cite{GH83} allows us to bypass deep neural network performing synchronized dimensional reduction and SINDy approach in \cite{Brunton19}. It is shown in Table \ref{table:table-result} that performing dimensional reduction and SINDy independently does not result in any significant reconstruction loss. Adding priors rooted in the theory of dynamical systems can improve and distinguish the SINDy algorithm from ordinary machine learning/feature engineering algorithms. Another prior that can be explored to further improve the quality of the discovered SINDy equations to reproduce \eqref{normalform} more faithfully is to modify the LASSO algorithm such that solves \eqref{eq:sindy-cost} while prioritizing the discovery of the linear coefficients, followed by the cubic coefficients, and lastly the remaining (quadratic) terms in the library. 

\begin{acknowledgements}
The authors acknowledge partial support for this work from Natural Sciences and Engineering Research Council~(NSERC) Discovery grant 50503-10802, TECSIS /Fields-CQAM Laboratory for Inference and Prediction, and NSERC-CRD grant 543433-19. The authors are also grateful to Mr. Yiming Meng his contribution in the viscous MG simulation development.
\end{acknowledgements}

\section*{Conflict of interest}
The authors declare that they have no conflict of interest.

\newpage

\appendix
\section{Appendix: Other Discovered Equations}
\label{sec:appendix}
\subsection{PCA and SINDy}
\begin{figure}[htb]
    \centering
    \includegraphics[width=3.5in]{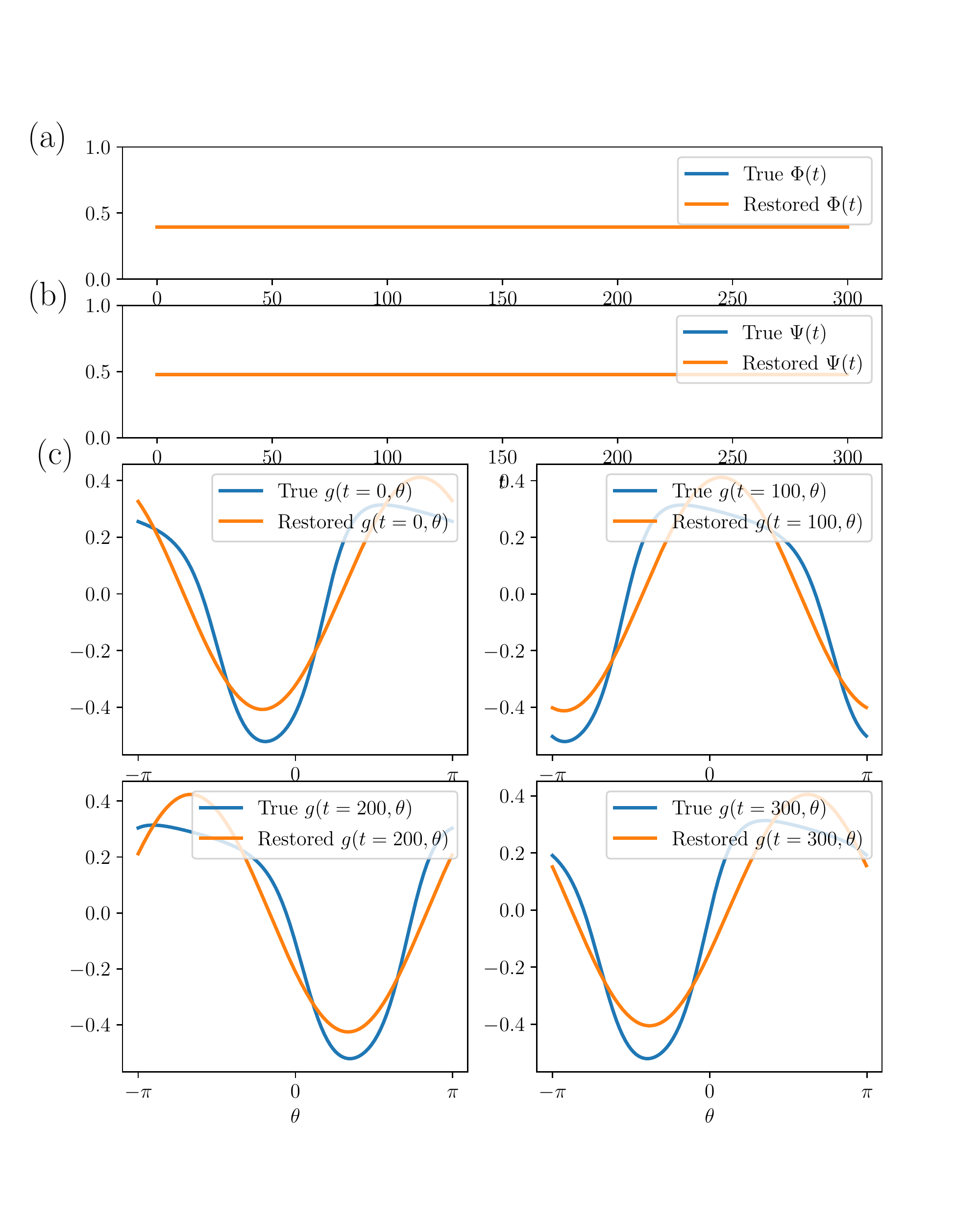}
    \caption{Reconstruction of viscous MG stall dynamics from multiplying ${\bf P}_{m=2}$ with numerical integration result of \eqref{pod-sindy}. (a) $\Phi(t)$ reconstruction result. (b) $\Psi(t)$ reconstruction result. (c) $g(t,\theta)$ reconstruction results at $t=0, 100, 200, 300$.}
    \label{pod-reconstruction-app}
\end{figure}

The SINDy equations that capture the symmetry of the cubic terms are obtained using a LASSO threshold of  $\alpha = 0.11$, which outputs a system of ODEs with 8 coefficients and test $\||{\bf \dot{X}}\||$ score of 0.9999
\begin{eqnarray}
\dot{x}_1&=& -0.000144 x_1^3  -0.008137 x_1^2 x_2  -0.000144 x_1 x_2^2 -0.008136 x_2^3
\nonumber \\
\dot{x}_2 &=& 0.008136 x_1^3 -0.000144 x_1^2 x_2 + 0.008136 x_1 x_2^2 -0.000144 x_2^3.
\label{pod-sindy-app}
\end{eqnarray}
The normal form coefficients $(-0.000144, 0.008136)$ are visibly detected. Due to the the higher $\alpha$ value, the linear terms are not captured. Reconstruction result for a chosen random dataset is shown in Figure \ref{pod-reconstruction-app}.

\subsection{Regularized Linear Autoencoder and SINDy}
\begin{figure}[htb]
    \centering
    \includegraphics[width=3.5in]{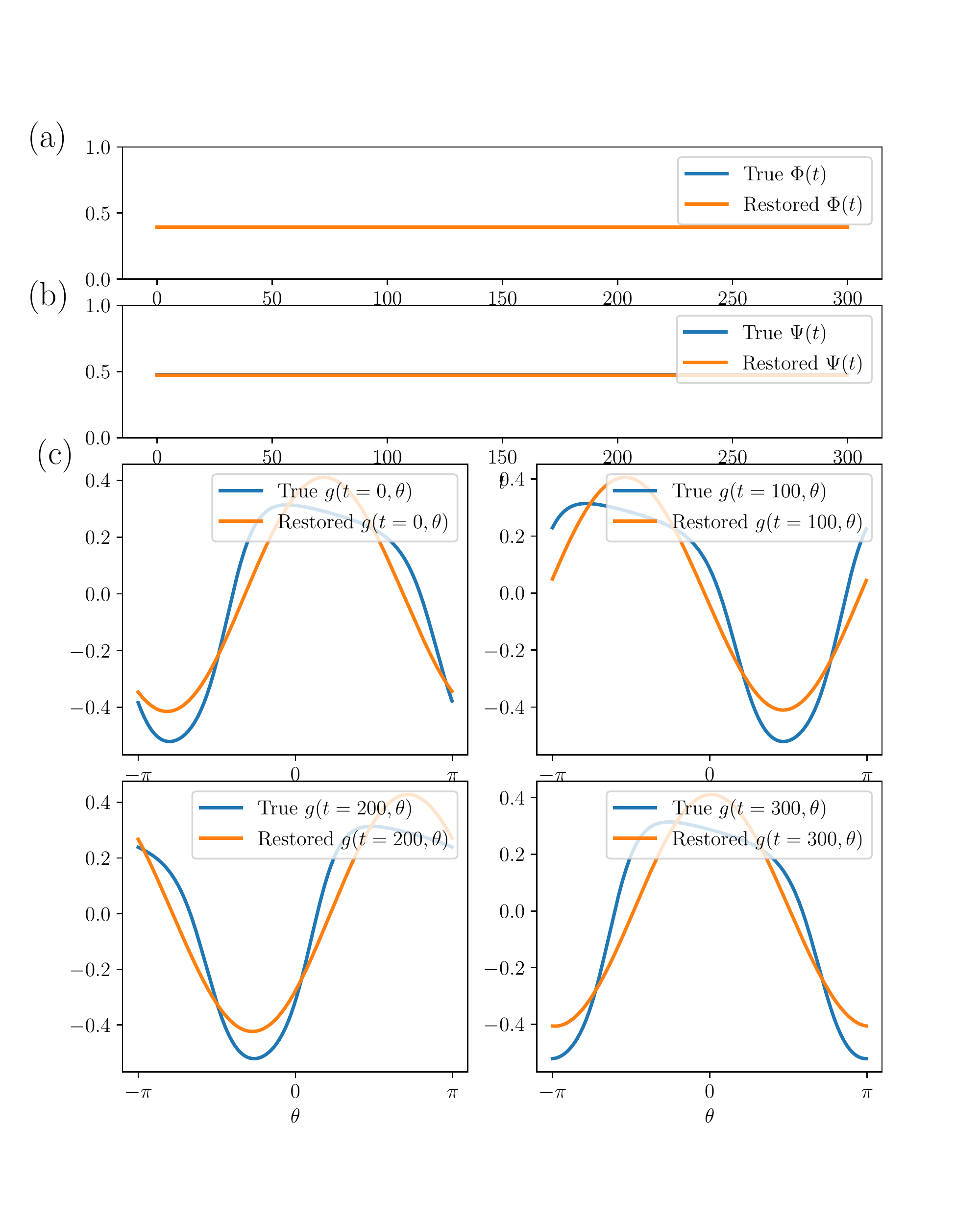}
    \caption{Reconstruction of viscous MG stall dynamics from feeding numerical integration result of \eqref{svdnn-sindy} into the decoder part of a regularized linear autoencoder. (a) $\Phi(t)$ reconstruction result. (b) $\Psi(t)$ reconstruction result. (c) $g(t,\theta)$ reconstruction results at $t=0, 100, 200, 300$.}
    \label{svdnn-reconstruction-app}
\end{figure}

The SINDy equations that capture the symmetry of the cubic terms are obtained using a LASSO threshold of  $\alpha = 0.30$. The output is a system of ODEs with 9 coefficients and test $\||{\bf \dot{X}}\||$ score of 0.9999
\begin{eqnarray}
\dot{x}_1&=& -0.000440 x_2^2 -0.000145 x_1^3 + 0.008138 x_1^2 x_2 -0.000143 x_1 x_2^2
\nonumber \\
&& + 0.008138 x_2^3
\nonumber \\
\dot{x}_2&=& -0.008136 x_1^3 -0.000144 x_1^2 x_2 -0.008139 x_1 x_2^2 -0.000144 x_2^3.
\label{svdnn-sindy-app}
\end{eqnarray}
The normal form coefficients $(-0.000145, -0.008138)$ are also visibly detected, albeit with wider deviation in values compared to \eqref{pod-sindy} and an additional quadratic term in the first equation. Due to the the higher $\alpha$ value, the linear terms are not captured. Additionally, at this $\alpha$ value, the resulting coefficients found are more consistent even with the different optimal encoders, unlike \eqref{svdnn-sindy}. Reconstruction result for a chosen random dataset is shown in Figure \ref{svdnn-reconstruction-app}.

\subsection{NLPCA Autoencoder and SINDy}
\begin{figure}[htb]
    \centering
    \includegraphics[width=3.5in]{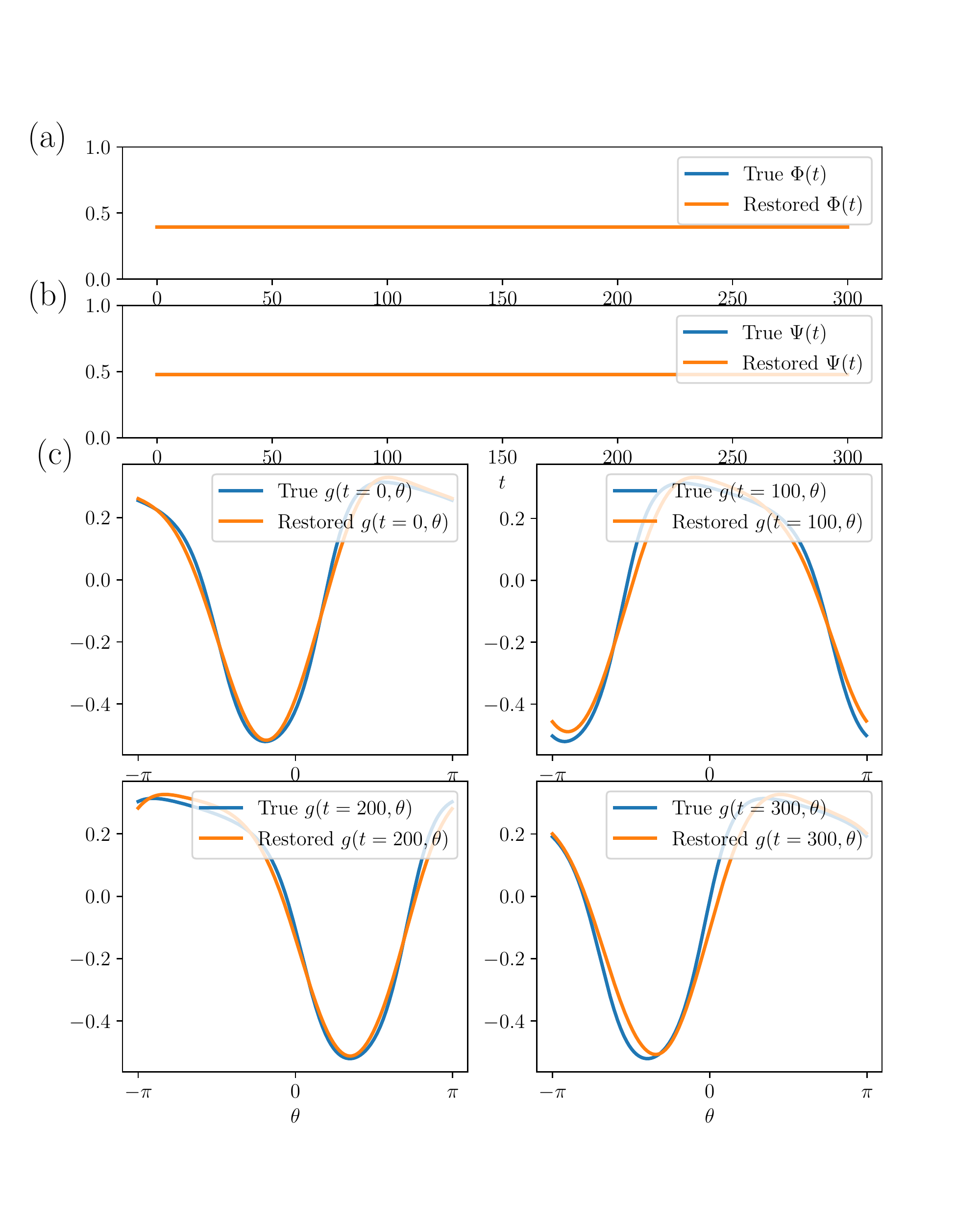}
    \caption{Reconstruction of viscous MG stall dynamics from feeding numerical integration result of \eqref{nlpca-sindy} into the decoder part of an NLPCA autoencoder. (a) $\Phi(t)$ reconstruction result. (b) $\Psi(t)$ reconstruction result. (c) $g(t,\theta)$ reconstruction results at $t=0, 100, 200, 300$.}
    \label{nlpca-reconstruction-app}
\end{figure}

The best SINDy regression is obtained using a LASSO threshold of  $\alpha = 0.60$. The output is a system of ODEs with 11 coefficients and test $\||{\bf \dot{X}}\||$ score of 0.9998 which does not satisfy the normal form. A representative equation (since the outcome is always random) is 
\begin{eqnarray}
\dot{x}_1&=& -0.001031 x_2^2  -0.000091 x_1^3 + 0.007704 x_1^2 x_2 + 0.000429 x_1 x_2^2
\nonumber \\
&& + 0.006863 x_2^3
\nonumber \\
\dot{x}_2&=& -0.001664 x_1^2  -0.000580 x_2^2  -0.007607 x_1^3 -0.000473 x_1^2 x_2
\nonumber \\
&&-0.007470 x_1 x_2^2 -0.000237 x_2^3.
\label{nlpca-sindy-app}
\end{eqnarray}
Reconstruction result for a chosen random dataset is shown in Figure \ref{nlpca-reconstruction-app}.

%

\end{document}